\documentclass[12pt]{article}

\usepackage[utf8]{inputenc}
\usepackage[a4paper, margin=0.8in]{geometry}
\usepackage{amsmath, amssymb, amsfonts}
\usepackage{graphicx}
\usepackage{float}
\usepackage{times}
\usepackage{hyperref}
\hypersetup{
  colorlinks=true,
  linkcolor=blue,
  citecolor=blue,
  urlcolor=blue
}
\usepackage{setspace}
\usepackage{subcaption}
\begin{document}

\title{\Large Characteristic‑Based Physics‑Informed Neural Networks for Advection Equations with a Focus on Discontinuous Solutions}
\author{
Omid Khosravi \and
Mehdi Tatari
}

\date{}

\maketitle

\begin{center}
Department of Mathematical Sciences, Isfahan University of Technology, Isfahan, 84156-83111, Iran
\end{center}

\begin{abstract}
This paper presents a characteristic-based loss function for physics-informed neural networks (PINNs) that circumvents automatic differentiation and significantly reduces training cost for advection equations. 
However, standard PINNs often fail to accurately resolve discontinuous solutions, as their global smoothness prior and spectral bias lead to smeared shocks or spurious oscillations near sharp gradients. 
To accurately capture these discontinuities, the density of sampling points in their vicinity must be sufficiently high, and an appropriate sampling strategy must be employed. 
For problems involving discontinuous initial and boundary conditions, several complementary techniques are introduced: 
an adaptive sampling strategy that concentrates collocation points along characteristic paths, 
a Fourier feature mapping with two-stage training and adaptive weighting to mitigate spectral bias, 
and an adaptive median filter applied to the spatial data at each time instant, which suppresses spurious oscillations while preserving sharp features such as corners and peaks. 
The effectiveness of the proposed framework is demonstrated through numerical experiments, showing improved accuracy and reduced computational effort in approximating discontinuous solutions.
\end{abstract}


\section{Introduction}

In recent decades, artificial neural networks (ANNs) have become a powerful class of computational models, capable of solving a wide variety of problems when their architecture and training procedure are appropriately designed.

Physics-informed neural networks (PINNs) \cite{raissi2019physics} combine machine learning techniques with physical laws to numerically solve partial differential equations (PDEs). In this framework, a parametric neural network \(\hat{u}(\mathbf{x},t;\boldsymbol{\theta})\) is employed as an approximation to the unknown solution \(u(\mathbf{x},t)\). The network parameters \(\boldsymbol{\theta}\) are trained such that the approximation simultaneously satisfies the governing equation together with the prescribed initial and boundary conditions. The general problem can be written as
\begin{equation*}
\mathcal{N}[u](\mathbf{x},t)=0,\qquad \mathbf{x}\in\Omega,\ t\in[0,T],
\end{equation*}
subject to
\begin{equation*}
u(\mathbf{x},0)=u_0(\mathbf{x}),\qquad \mathbf{x}\in\Omega,
\end{equation*}
and
\begin{equation*}
\mathcal{B}[u](\mathbf{x},t)=g(\mathbf{x},t),\qquad (\mathbf{x},t)\in\partial\Omega\times[0,T],
\end{equation*}
where \(\mathcal{N}\) denotes the differential operator of the governing equation, \(\mathcal{B}\) is the boundary operator, \(\Omega\subset\mathbb{R}^d\) is the spatial domain, and \(T\) is the final time.

The training process minimizes a composite loss function that enforces the PDE residual, the initial condition, and the boundary condition simultaneously. Denoting the collocation points for the PDE residual by \(\{\mathbf{x}^{i}_{\mathrm{PDE}},t^{i}_{\mathrm{PDE}}\}_{i=1}^{n_{\mathrm{PDE}}}\), the initial condition points by \(\{\mathbf{x}^{i}_{\mathrm{IC}}\}_{i=1}^{n_{\mathrm{IC}}}\), and the boundary points by \(\{\mathbf{x}^{i}_{\mathrm{BC}},t^{i}_{\mathrm{BC}}\}_{i=1}^{n_{\mathrm{BC}}}\), the individual loss components are defined as
\begin{align*}
\mathcal{L}_{\mathrm{PDE}}(\boldsymbol{\theta})
&= \frac{1}{n_{\mathrm{PDE}}}\sum_{i=1}^{n_{\mathrm{PDE}}}
\big|\mathcal{N}[\hat{u}](\mathbf{x}^{i}_{\mathrm{PDE}},t^{i}_{\mathrm{PDE}};\boldsymbol{\theta})\big|^{2},\\
\mathcal{L}_{\mathrm{IC}}(\boldsymbol{\theta})
&= \frac{1}{n_{\mathrm{IC}}}\sum_{i=1}^{n_{\mathrm{IC}}}
\big|\hat{u}(\mathbf{x}^{i}_{\mathrm{IC}},0;\boldsymbol{\theta})-u_{0}(\mathbf{x}^{i}_{\mathrm{IC}})\big|^{2},\\
\mathcal{L}_{\mathrm{BC}}(\boldsymbol{\theta})
&= \frac{1}{n_{\mathrm{BC}}}\sum_{i=1}^{n_{\mathrm{BC}}}
\big|\mathcal{B}[\hat{u}](\mathbf{x}^{i}_{\mathrm{BC}},t^{i}_{\mathrm{BC}};\boldsymbol{\theta})
- g(\mathbf{x}^{i}_{\mathrm{BC}},t^{i}_{\mathrm{BC}})\big|^{2}.
\end{align*}
The total loss is a weighted sum
\begin{equation*}
\mathcal{L}(\boldsymbol{\theta})
= \lambda_{\mathrm{PDE}}\,\mathcal{L}_{\mathrm{PDE}}(\boldsymbol{\theta})
+ \lambda_{\mathrm{IC}}\,\mathcal{L}_{\mathrm{IC}}(\boldsymbol{\theta})
+ \lambda_{\mathrm{BC}}\,\mathcal{L}_{\mathrm{BC}}(\boldsymbol{\theta}),
\end{equation*}
with non-negative weights \(\lambda_{\mathrm{PDE}}, \lambda_{\mathrm{IC}}, \lambda_{\mathrm{BC}}\). Automatic differentiation (AD) \cite{baydin2018automatic} is used to compute the spatial and temporal derivatives that appear in \(\mathcal{N}[\hat{u}]\), removing the need for an explicit mesh. While this approach is flexible and general, the AD-based residual can be computationally expensive and may amplify high‑frequency errors near sharp gradients or discontinuities.

Although PINNs are efficient and versatile, they encounter difficulties with solutions that exhibit strong nonlinearities, sharp gradients, or discontinuities. Several modifications have been proposed to improve their performance, such as domain decomposition strategies (DPINN \cite{dwivedi2021distributed}, cPINN \cite{jagtap2020conservative}, XPINN \cite{jagtap2020extended}). Another well-known limitation of neural networks is spectral bias \cite{rahaman2019spectral} — the tendency to learn low-frequency components first, which hinders the reconstruction of fine-scale features. Random Fourier features \cite{tancik2020fourier} have been introduced to mitigate this issue by mapping the input to a higher-dimensional sinusoidal space. A variant of this idea, with trainable frequencies and adaptive weighting, is employed in the present work to better capture discontinuities.

The advection equation is one of the simplest yet most challenging PDEs from a numerical point of view \cite{morton2005numerical}. In its most general form,
\begin{equation}\label{eq:advection}
\frac{\partial u}{\partial t} + \mathbf{a}(\mathbf{x},t,u)\cdot\nabla u = f(\mathbf{x},t), \qquad \mathbf{x}\in\Omega\subset\mathbb{R}^d,\; t\ge 0,
\end{equation}
where \(u(\mathbf{x},t)\) is the transported quantity, \(\mathbf{a}\) is the velocity field (which may depend on \(u\)), and \(f\) is a source term. The problem is linear if \(\mathbf{a}\) is independent of \(u\). The exact solution is constant along characteristic curves, and numerical methods such as upwind and Lax--Wendroff have been developed to approximate it. The upwind scheme preserves the maximum principle under a CFL condition, while Lax--Wendroff does not always do so \cite{hosseini2019splitting}. The presence of discontinuities in the initial or boundary conditions makes the problem significantly more difficult, as the jumps propagate along characteristics.

In this paper, we propose to exploit the method of characteristics directly inside the PINN loss function, thereby avoiding automatic differentiation for the advection operator and substantially reducing the training cost. This characteristic-based loss is general and suitable for both smooth and discontinuous problems. To further improve the treatment of discontinuities, we introduce a suite of complementary techniques: an adaptive sampling strategy that places collocation points along characteristic paths near the jump, a trainable Fourier feature mapping with a two‑stage training procedure and adaptive loss weighting, and a spatial median filter that eliminates spurious oscillations while preserving discontinuities.

The remainder of the paper is organized as follows. Section~\ref{sec:char_loss} presents the characteristic-based loss function for linear and nonlinear advection. Section~\ref{sec:discont} describes the strategies tailored for discontinuous initial and boundary conditions: characteristic‑based sampling, Fourier feature mapping with adaptive weighting, and median filtering. Section~\ref{sec:experiments} reports numerical experiments, and Section~\ref{sec:conclusion} concludes the paper.

\section{Characteristic-Based Loss Function}
\label{sec:char_loss}

A fundamental property of the advection equation~\eqref{eq:advection} is that its solution evolves along characteristic curves. For a linear advection equation
\begin{equation}\label{eq:linear_adv}
\frac{\partial u}{\partial t} + \mathbf{a}(\mathbf{x},t)\cdot\nabla u = f(\mathbf{x},t),
\end{equation}
where \(\mathbf{a}\) depends only on \(\mathbf{x}\) and \(t\), the characteristics are determined by the system of ordinary differential equations
\begin{equation}\label{eq:char_ode}
\begin{cases}
\dfrac{d\mathbf{x}}{dt} = \mathbf{a}(\mathbf{x},t), \\[10pt]
\dfrac{du}{dt} = f(\mathbf{x},t).
\end{cases}
\end{equation}
For a given collocation point \((\mathbf{x}_{\mathrm{PDE}}^{i}, t_{\mathrm{PDE}}^{i})\), we trace the characteristic forward in time by solving \(\frac{d\mathbf{x}}{dt} = \mathbf{a}(\mathbf{x},t)\) numerically. Because the velocity field is independent of \(u\), the characteristic paths can be computed once, before training, with high accuracy using methods such as the classical fourth-order Runge--Kutta (RK4) or the implicit Radau method. Using a half-step \(h/2\) and two consecutive steps yields both the midpoint \(\overline{\mathbf{x}}_{\mathrm{PDE}}^{i}\) at time \(t_{\mathrm{PDE}}^{i}+h/2\) and the characteristic foot \(\tilde{\mathbf{x}}_{\mathrm{PDE}}^{i}\) at time \(t_{\mathrm{PDE}}^{i}+h\). However, for greater accuracy, one may compute these two positions using a smaller step size.

Precomputing these points completely eliminates spatial automatic differentiation during training and substantially reduces the per‑iteration cost.

Integrating the second ODE in~\eqref{eq:char_ode} along the characteristic and applying the midpoint rule yields the approximate relation
\begin{equation}\label{eq:char_relation}
\frac{u(\tilde{\mathbf{x}}_{\mathrm{PDE}}^{i}, t_{\mathrm{PDE}}^{i}+h) - u(\mathbf{x}_{\mathrm{PDE}}^{i}, t_{\mathrm{PDE}}^{i})}{h} \approx f\!\left(\overline{\mathbf{x}}_{\mathrm{PDE}}^{i},\, t_{\mathrm{PDE}}^{i}+\tfrac{h}{2}\right).
\end{equation}
If \(\tilde{\mathbf{x}}_{\mathrm{PDE}}^{i}\) and \(\overline{\mathbf{x}}_{\mathrm{PDE}}^{i}\) are computed using a method of at least second-order accuracy, the above approximation is \(\mathcal{O}(h^{2})\).  
Moreover, when the source term vanishes (\(f = 0\)), the accuracy of the characteristic relation reduces to the accuracy of the computed foot position \(\tilde{\mathbf{x}}_{\mathrm{PDE}}^{i}\), as the relation simplifies to \(u(\tilde{\mathbf{x}}, t+h) = u(\mathbf{x}, t)\).

We therefore define the \emph{characteristic-based PDE loss} as
\begin{equation}\label{eq:char_loss}
\mathcal{L}_{\mathrm{PDE}}(\boldsymbol{\theta})
= \frac{1}{n_{\mathrm{PDE}}}\sum_{i=1}^{n_{\mathrm{PDE}}}
\Bigg|
\frac{\hat{u}\big(\tilde{\mathbf{x}}_{\mathrm{PDE}}^{i}, t_{\mathrm{PDE}}^{i}+h; \boldsymbol{\theta}\big)
- \hat{u}\big(\mathbf{x}_{\mathrm{PDE}}^{i}, t_{\mathrm{PDE}}^{i}; \boldsymbol{\theta}\big)}{h}
- f\big(\overline{\mathbf{x}}_{\mathrm{PDE}}^{i}, t_{\mathrm{PDE}}^{i} + \tfrac{h}{2}\big)
\Bigg|^{2}.
\end{equation}

This loss function completely bypasses automatic differentiation. Evaluating it requires only two forward passes through the network (original points and characteristic feet), whereas the standard PINN residual requires the construction of the derivative graph and an additional pass through it, which is significantly heavier. The computational cost is therefore substantially lower. Moreover, since spatial derivatives are never computed, high-frequency errors near discontinuities are not amplified.

The step size \(h\) is a hyperparameter. Because neural network training is typically performed in float32 arithmetic, \(h\) in the range \(10^{-3}\) to \(10^{-2}\) is recommended to balance truncation and round-off errors.

\medskip
\noindent\textbf{Extension to nonlinear advection.}
When the velocity field depends on \(u\), i.e., \(\mathbf{a} = \mathbf{a}(\mathbf{x},t,u)\), the characteristics are no longer independent of the solution and cannot be precomputed. In this case we propose to solve the characteristic system
\[
\frac{d\mathbf{x}}{dt} = \mathbf{a}(\mathbf{x},t,u), \qquad \frac{du}{dt} = f(\mathbf{x},t)
\]
using the explicit Euler method. Starting from \((\mathbf{x}_{\mathrm{PDE}}^{i}, t_{\mathrm{PDE}}^{i})\) and using the current network prediction \(\hat{u}(\mathbf{x}_{\mathrm{PDE}}^{i}, t_{\mathrm{PDE}}^{i};\boldsymbol{\theta})\), a single forward Euler step of size \(h\) is taken:
\[
\tilde{\mathbf{x}}_{\mathrm{PDE}}^{i} = \mathbf{x}_{\mathrm{PDE}}^{i} + h\,\mathbf{a}\big(\mathbf{x}_{\mathrm{PDE}}^{i}, t_{\mathrm{PDE}}^{i}, \hat{u}(\mathbf{x}_{\mathrm{PDE}}^{i}, t_{\mathrm{PDE}}^{i};\boldsymbol{\theta})\big).
\]
The corresponding discrete relation is
\[
\frac{u(\tilde{\mathbf{x}}_{\mathrm{PDE}}^{i}, t_{\mathrm{PDE}}^{i}+h) - u(\mathbf{x}_{\mathrm{PDE}}^{i}, t_{\mathrm{PDE}}^{i})}{h} \approx f(\mathbf{x}_{\mathrm{PDE}}^{i}, t_{\mathrm{PDE}}^{i}),
\]
and the characteristic-based loss becomes
\begin{equation}\label{eq:char_loss_nonlinear}
\mathcal{L}_{\mathrm{PDE}}(\boldsymbol{\theta})
= \frac{1}{n_{\mathrm{PDE}}}\sum_{i=1}^{n_{\mathrm{PDE}}}
\Bigg|
\frac{\hat{u}\big(\tilde{\mathbf{x}}_{\mathrm{PDE}}^{i}, t_{\mathrm{PDE}}^{i}+h; \boldsymbol{\theta}\big)
- \hat{u}\big(\mathbf{x}_{\mathrm{PDE}}^{i}, t_{\mathrm{PDE}}^{i}; \boldsymbol{\theta}\big)}{h}
- f(\mathbf{x}_{\mathrm{PDE}}^{i}, t_{\mathrm{PDE}}^{i})
\Bigg|^{2}.
\end{equation}
Because Euler's method is only first-order accurate, this approximation is \(\mathcal{O}(h)\). The dependence of \(\tilde{\mathbf{x}}_{\mathrm{PDE}}^{i}\) on \(\hat{u}\) makes the computational graph roughly twice as long as a direct evaluation; however, it remains significantly lighter than the graph required for automatic differentiation.

The characteristic-based loss is therefore a general and efficient alternative to the standard PINN residual, suitable for both smooth and discontinuous advection problems.

Moreover, if Neumann boundary conditions are prescribed, the normal derivative can be approximated by a forward difference in the direction opposite to the outward normal, with arbitrary order of accuracy (e.g., \((u(\mathbf{x}) - u(\mathbf{x}-h\mathbf{n}))/h\), which is first‑order). This requires only additional network evaluations at shifted interior points and avoids automatic differentiation entirely.
\section{Strategies for Handling Discontinuous Initial and Boundary Conditions}
\label{sec:discont}

In this section we describe three complementary techniques that improve the PINN approximation when the initial or boundary data contain jump discontinuities. These strategies address the placement of collocation points, the representation of high‑frequency features, and the suppression of spurious oscillations.

\subsection{Characteristic-Based Sampling}
\label{sec:adaptive_sampling}

Neural networks have difficulty learning sharp, discontinuous behaviors because they inherently generate smooth functions; they approximate a jump by a steep continuous transition, which introduces significant error and slows convergence. For problems with discontinuous conditions, it is recommended to place the corresponding training points at a small distance from the exact location of the jump, avoiding evaluation where the function is not defined or the gradient would be infinite. This practice controls the slope of the approximated transition and improves training.

Furthermore, maintaining a suitable density of interior points around the propagation path of the discontinuity is crucial. For linear equations, this can be achieved by exploiting the characteristic structure. Points are selected in the vicinity of the discontinuity on the initial and boundary data, and the ODE
\[
\frac{d\mathbf{x}}{dt} = \mathbf{a}(\mathbf{x},t)
\]
is solved with a suitable numerical method over the time interval. Tracking the discontinuity along the characteristics yields a set of interior collocation points that lie precisely along the jump path. To avoid over-concentration on a single curve, to improve generalization, and to ensure a small offset from the exact discontinuity (which helps the network learn the steep transition), a small amount of Gaussian noise is added to the space–time coordinates of these points. This strategy confines high point density to the immediate neighborhood of the discontinuity, preventing an unnecessary increase in the total number of collocation points and significantly reducing training cost.

In the nonlinear case, the same idea applies, but one must solve the full characteristic system
\[
\frac{d\mathbf{x}}{dt} = \mathbf{a}(\mathbf{x},t,u), \qquad \frac{du}{dt} = f(\mathbf{x},t)
\]
to track the discontinuity. Compared to the linear case, this process may involve greater cost and additional challenges.

Even when the solution is globally continuous but the initial condition contains a sharp peak (e.g., a narrow Gaussian), the characteristic‑based sampling strategy remains effective: the peak location can be transported along the characteristics, providing a set of collocation points that closely follow the moving high‑gradient region.
\subsection{Fourier Feature Mapping}
\label{sec:fourier}

Random Fourier features enhance the ability of neural networks to model rapidly oscillatory functions, local discontinuities, and high-frequency components, counteracting the spectral bias of multilayer perceptrons (the tendency to learn low frequencies first \cite{rahaman2019spectral}).

An input vector \(\mathbf{m}\in\mathbb{R}^d\) is mapped nonlinearly before being passed to the network:
\begin{equation}\label{eq:fourier}
\boldsymbol{\gamma}(\mathbf{m}) = \big[
\cos(\mathbf{B}\mathbf{m}),\;\sin(\mathbf{B}\mathbf{m})
\big]^\top \in\mathbb{R}^{2D},
\end{equation}
where \(D\) is the number of sine–cosine pairs and \(\mathbf{B}\in\mathbb{R}^{D\times d}\) encodes the frequency scales and directions. The rows of \(\mathbf{B}\) are typically sampled from a Gaussian distribution \(\mathcal{N}(0,\sigma^2)\). For a time-dependent problem in \(d\) spatial dimensions, \(\mathbf{m}=[\mathbf{x},\;t]^\top \in \mathbb{R}^{d+1}\). The parameter \(\sigma\) controls the frequency spectrum: larger values improve the capture of fine scales, but excessively large values may cause instability and Gibbs-like artifacts.

In addition to improving representational capacity, Fourier features often accelerate optimizer convergence. The transformed input exposes both low- and high-frequency components in a quasi-linear fashion, enabling the model to represent localized structures without increasing network depth.

\subsubsection{Trainable Fourier Feature Mapping}
\label{sec:trainable_fourier}

Training the entries of \(\mathbf{B}\) before training the network produces a more targeted feature transform that is particularly beneficial for discontinuities. We propose a two-stage scheme: first, the Fourier mapping parameters \(\boldsymbol{\theta}_1\) are optimized while the network parameters \(\boldsymbol{\theta}_2\) are frozen; then \(\boldsymbol{\theta}_1\) is frozen and \(\boldsymbol{\theta}_2\) is trained. The model is written as \(\hat u(\mathbf{x},t;\boldsymbol{\theta}) = N(\gamma(\mathbf{x},t;\boldsymbol{\theta}_1);\boldsymbol{\theta}_2)\). This procedure significantly accelerates training in the presence of discontinuities and helps rebalance the loss components. The network weights are initialized using the Glorot \cite{glorot2010understanding} scheme, and biases are set to zero. Longer training in the first stage tends to balance the discontinuous components relative to the smooth ones, expediting overall convergence. Once the Fourier mapping has stabilized, the network parameters are trained separately; the training data remain identical across both stages.

\subsubsection{Adaptive Fourier Weighting}
\label{sec:adaptive_weighting}

Since the primary source of imbalance in the loss function is the discontinuity itself, we assign a significantly larger weight to the loss component associated with the discontinuous part (initial or boundary conditions) during the first stage — a strategy we call \emph{adaptive Fourier weighting}. This focuses the mapping optimization on regions containing jumps, leading to faster and more balanced learning of these features in the second stage. The resulting feature space is better aligned with the requirements of modeling discontinuous solutions. Moreover, adaptive weighting reduces the likelihood of convergence to undesirable local minima and empirically improves training speed significantly, while alleviating the need for highly complex architectures.

In the second stage (network‑parameter optimization), alternative weighting schemes such as gradient-norm‑based methods or NTK‑based weighting can be employed; nevertheless, it is advisable that the loss component for the discontinuous part still receive a larger weight to maintain balance. The combination of two‑stage training and adaptive weighting effectively balances the loss components and substantially reduces overall training time.

\subsection{Median Filtering}
\label{sec:median}

The strategies above can significantly accelerate training, but the combination of Fourier features and data scarcity near discontinuities may still lead to overfitting and unwanted oscillations. To mitigate these effects, we apply a median filter.

The median filter is a non‑linear smoothing operator that replaces each value with the median of its neighbors. In \(d\) dimensions, for an array entry with coordinates \(\mathbf{i}\), a hypercube (or hyperrectangle) of odd size \(k\) is centered at \(\mathbf{i}\). The filtered value is
\[
u^{\text{filtered}}(\mathbf{i}) = \mathcal{M}\Big(\big\{ u(\mathbf{j}) \mid \mathbf{j} \in \mathcal{N}(\mathbf{i}) \big\}\Big),
\]
where \(\mathcal{N}(\mathbf{i})\) is the set of indices within the window and \(\mathcal{M}\) is the median operator.

When solving the advection equation, we apply the median filter exclusively to the spatial data at each fixed time instant. At every time step, the spatial field \(\hat{u}(\mathbf{x},t;\boldsymbol{\theta})\) is filtered once. The window size \(k\) is a hyperparameter controlling the degree of smoothing. To avoid boundary artefacts, points whose window extends outside the domain are left unfiltered. The median filter preserves edges and discontinuities while effectively removing impulse-like oscillations.

It should be noted that if the exact solution itself contains point‑like peaks or corner discontinuities, median filtering can degrade these physical features. To address this, we exploit the same characteristic‑tracking idea used in the adaptive sampling strategy (Section~\ref{sec:adaptive_sampling}): the initial corners of the discontinuity or sharp peak points are transported along the characteristics. At each time step, a small neighborhood around each of these tracked points is marked, and the median filter is not applied inside those regions. In this way, sharp corners and peaks are preserved while the rest of the domain benefits from the denoising effect of the filter.

In many problems, especially where oscillations are artificial and caused by network error, this combined filtering strategy yields a significant improvement in approximation quality.

\section{Numerical Experiments}
\label{sec:experiments}

All experiments were run on Google Colab equipped with an NVIDIA Tesla T4 GPU, using PyTorch~2.11.0+cu128. For every comparative run, the neural network weights were initialized identically between the standard and characteristic‑based loss models to guarantee a fair comparison.

We present four test cases that systematically evaluate the proposed framework. The first is a two‑dimensional rigid‑body rotation with a discontinuous square pulse; it illustrates the complete set of techniques—characteristic‑based loss, adaptive sampling, Fourier features, and the adaptive median filter—and serves as the main benchmark for accuracy and computational cost. The second is a one‑dimensional advection of multiple blocks with a reduced network and limited data; it isolates the effects of the median filter, trainable Fourier features, and adaptive Fourier weighting in a minimal setting. The third example is a three‑dimensional smooth Gaussian rotation that contains no discontinuities; it provides a clean comparison of the two loss formulations without the confounding influence of jump‑specific techniques. The fourth example is a three‑dimensional nonlinear advection with a manufactured smooth solution; it demonstrates that the advantage of the characteristic‑based loss extends to nonlinear problems where the characteristic feet must be updated at each iteration.
\subsection{Two‑Dimensional Rigid‑Body Rotation}

Consider the linear advection equation
\[
\frac{\partial u}{\partial t} - (y-1)\frac{\partial u}{\partial x} + (x-1)\frac{\partial u}{\partial y} = 0,
\qquad (x,y) \in [0,2]^2,\; t \in [0,1],
\]
which describes a rigid rotation about the center \((1,1)\). The initial condition is a square pulse of unit value inside \([0.7,1.3]\times[0.7,1.3]\) and zero elsewhere; the boundary condition is homogeneous on the whole boundary. The exact solution is simply the initial square rotated rigidly, providing a clean visual and quantitative reference.

In this example, the neural network is a multilayer perceptron with four hidden layers of 20 neurons each and the hyperbolic tangent activation. The input is preprocessed by a Fourier feature mapping with \(D = 60\) sine–cosine pairs and \(\sigma = 3.0\); the mapping parameters are first optimized with Adam (learning rate \(2\times10^{-3}\), 4000 epochs) while the network weights are frozen. The learning rate in all Adam phases follows an exponential decay schedule, halving every 2000 epochs. In the second stage, the network weights are trained with Adam (8000 epochs). Finally, the model is fine‑tuned with the L‑BFGS optimizer (up to 3000 iterations) on a large precomputed batch. The characteristic‑based loss uses a step size \(h = 10^{-2}\). The unnormalized loss weights are set to \(\lambda_{\mathrm{IC}} = 10\), \(\lambda_{\mathrm{BC}} = 1\) and \(\lambda_{\mathrm{PDE}} = 1\). The median filter, applied on a uniform \(400\times400\) spatial grid at each time step, uses a window of size \(k = 7\); the corner‑preserving variant masks a circular neighborhood of radius \(0.02\) around each advected corner. For comparison, a standard PINN with the same network architecture is trained using the conventional PDE residual based on automatic differentiation; it follows the same optimization schedule and uses the same total number of collocation points as the final phase of the proposed method.

In this example, the Adam phase uses 4000 random points (resampled each epoch from a pool of 30 000) and 1000 characteristic points (drawn from the precomputed set of 10 100 points that track the discontinuity). The L‑BFGS phase employs 10 000 random and 10 000 characteristic points, all with precomputed characteristic feet. The standard PINN baseline uses a single batch of 20 000 uniformly distributed collocation points.

\begin{figure}[htbp]
    \centering
    \includegraphics[width=0.8\textwidth]{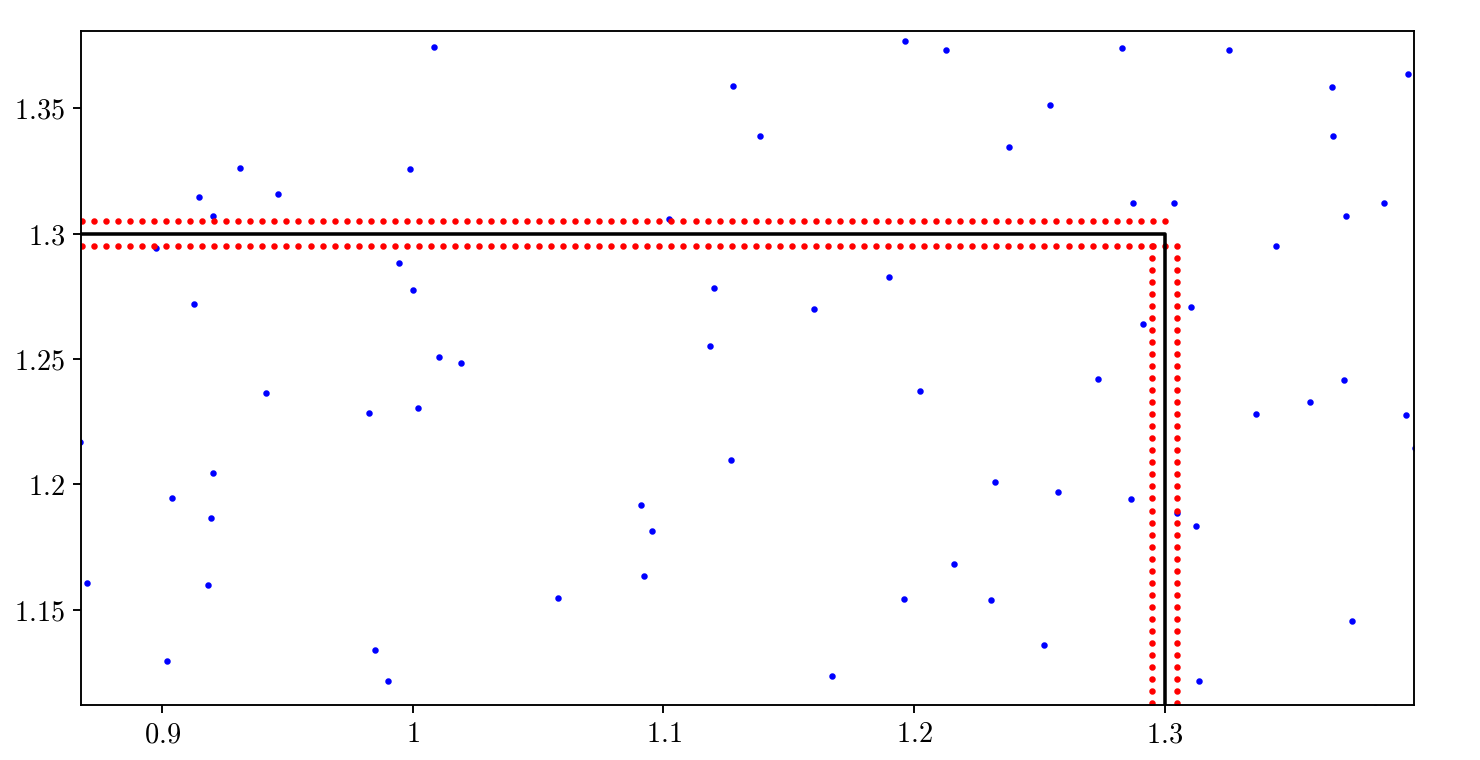}
    \caption{Initial condition training points. Blue: random points away from the jump. Red: points at distance \(\pm\delta\) from the discontinuity. The black line marks the true interface.}
    \label{fig:initial_points}
\end{figure}

To construct the initial condition training set, we place \(2000\) random points uniformly in the domain, excluding a narrow band of width \(2\delta = 0.01\) around the discontinuity, where the function is not differentiable. Additionally, \(500\) points are placed exactly at distance \(\delta = 0.005\) on both sides of the square boundary, following the local normal direction. Increasing this distance further improves the convergence speed and stability of the solution; however, an excessive increase may produce a non‑physical solution. The resulting point cloud is shown in Fig.~\ref{fig:initial_points}.

\begin{figure}[H]
    \centering
    \includegraphics[width=0.9\textwidth]{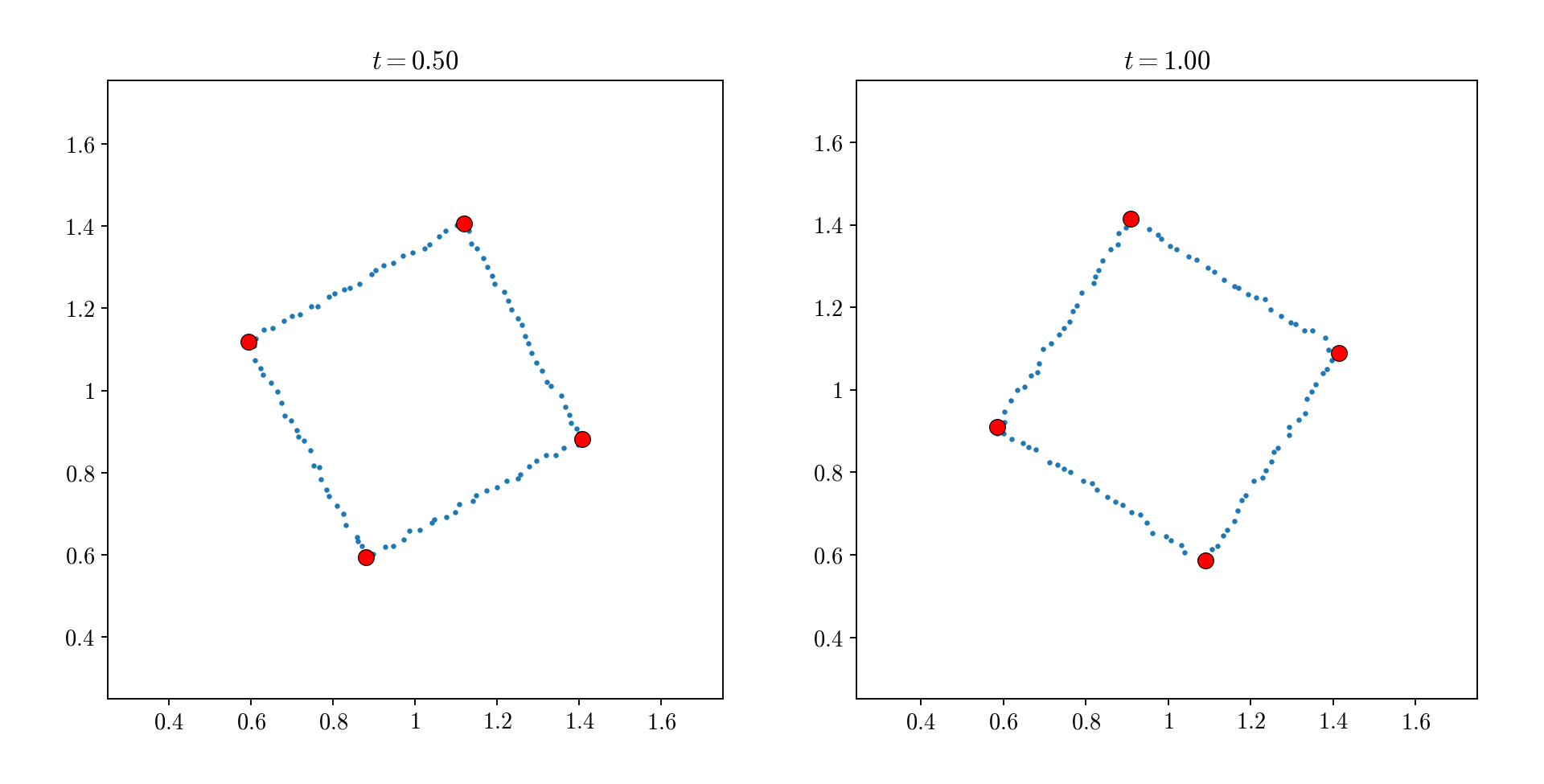}
    \caption{Characteristic‑based sampling points at times \(t = 0.5\) and \(t = 1.0\). Blue dots: points on the discontinuity path. Red circles: advected corners (radius \(0.02\)).}
    \label{fig:characteristic_points}
\end{figure}

To obtain interior collocation points that follow the discontinuity, we placed \(100\) points on the initial square boundary and solved the characteristic ODE \(\mathrm{d}\mathbf{x}/\mathrm{d}t = \mathbf{a}(\mathbf{x},t)\) using the Radau method. Due to the simplicity of the problem, the explicit Euler method would also provide sufficient accuracy. To avoid over‑concentration on a single curve, to improve generalization, and to ensure a small offset from the exact discontinuity (which helps the network learn the steep transition), a small amount of Gaussian noise (standard deviation \(0.005\)) is added to the space–time coordinates of the resulting points. The final set of points, shown in Fig.~\ref{fig:characteristic_points} (blue dots) for two selected time instants, lies near the moving discontinuity surface. The four corners of the square are also tracked along their own characteristics and are displayed as red circles; at each time step, a circular neighborhood of radius \(0.02\) around each corner is excluded from the median filter to preserve sharp geometric features. Without this adaptive sampling strategy, achieving a comparable level of accuracy would require a significantly larger number of uniformly distributed collocation points, leading to a substantial increase in training cost.

\begin{figure}[H]
    \centering
    \includegraphics[width=0.8\textwidth]{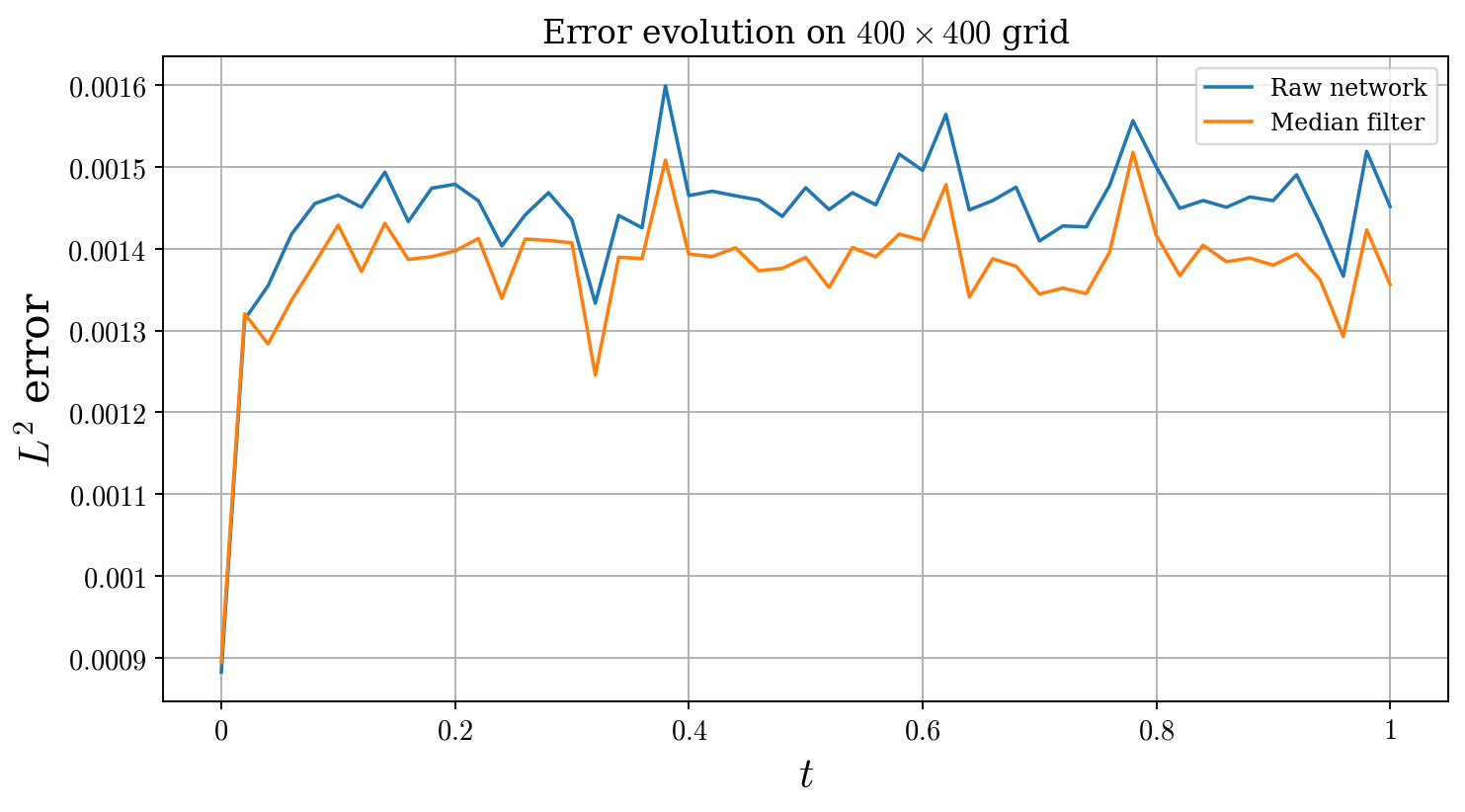}
    \caption{Squared \(L^2\) error evolution on a \(400\times400\) grid.}
    \label{fig:error_evolution_L2}
\end{figure}

\begin{figure}[H]
    \centering
    \includegraphics[width=0.8\textwidth]{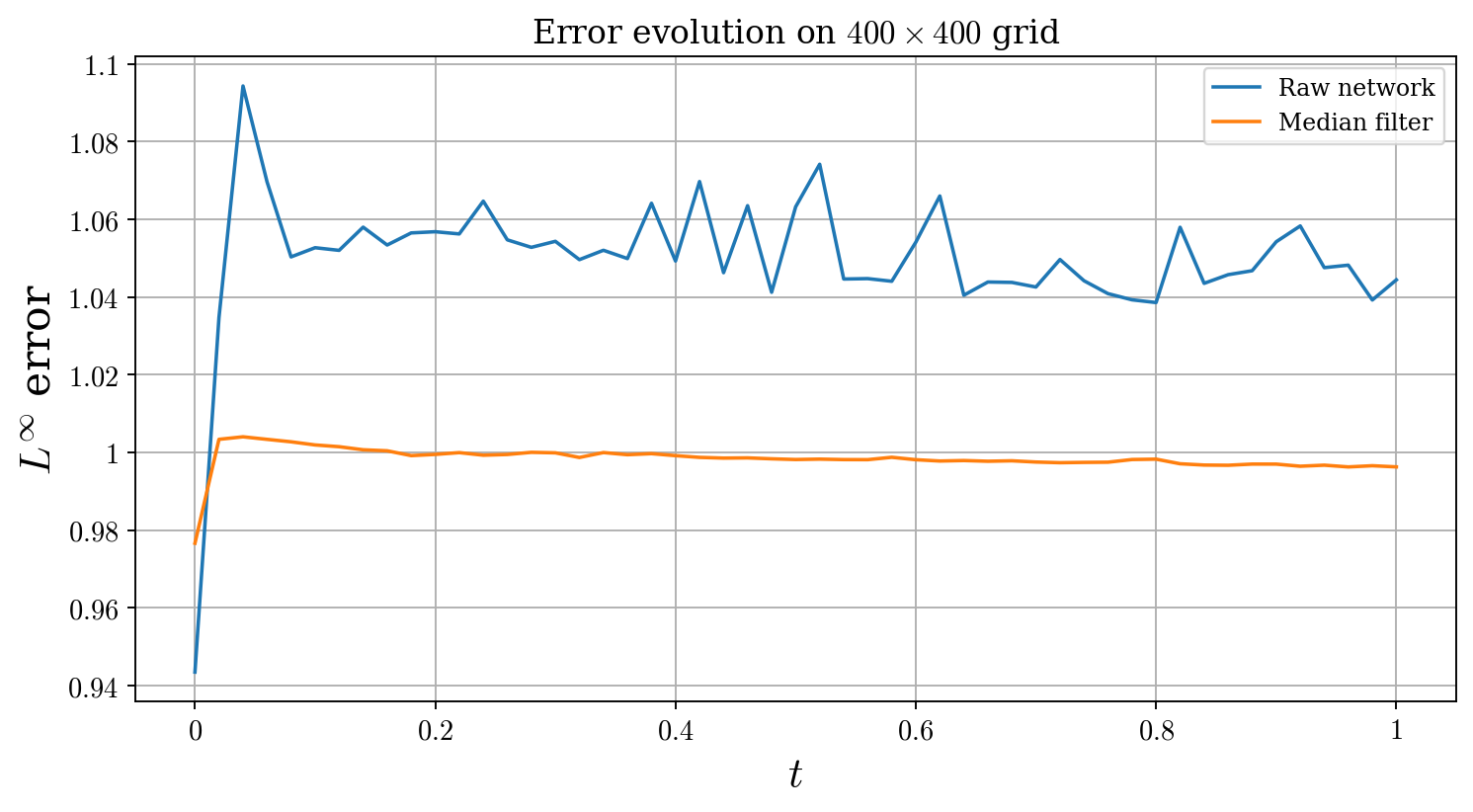}
    \caption{\(L^\infty\) error evolution on a \(400\times400\) grid.}
    \label{fig:error_evolution_Linf}
\end{figure}

Figure~\ref{fig:error_evolution_L2} shows the evolution of the squared \(L^2\) error on a \(400\times400\) uniform grid over time, for the raw network output, after a regular median filter, and after the corner‑preserving median filter. The corner‑preserving variant consistently yields the lowest error while preserving sharp features. A similar trend is observed for the \(L^\infty\) error in Fig.~\ref{fig:error_evolution_Linf}.

To quantify the variability due to random initialization, we repeated the experiment with three different random seeds. Table~\ref{tab:full_comparison} reports the per‑iteration Adam time, the total training time, and the squared \(L^2\) error (after the corner‑preserving median filter) for both the proposed characteristic‑based loss and the standard PINN loss. In this comparison, the only difference between the two approaches is the form of the PDE loss; the adaptive sampling strategy, the placement of initial and boundary points, and the corner‑preserving median filter are identical for both. The characteristic‑based loss significantly reduces the Adam iteration time; for the L‑BFGS optimizer (with a history size of \(500\) in our implementation) the relative speed‑up is smaller because the line‑search procedure dominates the per‑iteration cost, yet the proposed loss still retains an advantage. It is also worth noting that the gap in total training time becomes even more pronounced when a larger number of Adam iterations is employed.

\begin{table}[H]
  \centering
  \setlength{\tabcolsep}{12pt}
  \begin{tabular}{@{}ccccccc@{}}
    \hline
    Run & \multicolumn{2}{c}{Iter. time (s)} & \multicolumn{2}{c}{Total time (min)} & \multicolumn{2}{c}{Sq.\ \(L^2\) error} \\
    \cline{2-3} \cline{4-5} \cline{6-7}
    & Char. & Std. & Char. & Std. & Char. & Std. \\
    \hline
    1 & 0.0052 &  0.0092 & 180 & 253 & 1.46e-03 & 1.70e-03 \\
    2 & 0.0051 & 0.0093 & 177 & 250 & 1.57e-03 & 1.89e-03 \\
    3 & 0.0052 & 0.0091 & 182 & 247 & 9.73e-04 & 2.14e-3 \\
    \hline
  \end{tabular}
  \caption{Performance metrics over three random seeds. The Adam iteration time (seconds), total training time (minutes), and squared \(L^2\) error (on a \(400\times400\) grid at \(t=1\) after corner‑preserving median filter) are reported for the proposed characteristic‑based loss (Char.) and the standard PINN loss (Std.).}
  \label{tab:full_comparison}
\end{table}

\begin{figure}[H]
    \centering
    \includegraphics[width=\textwidth]{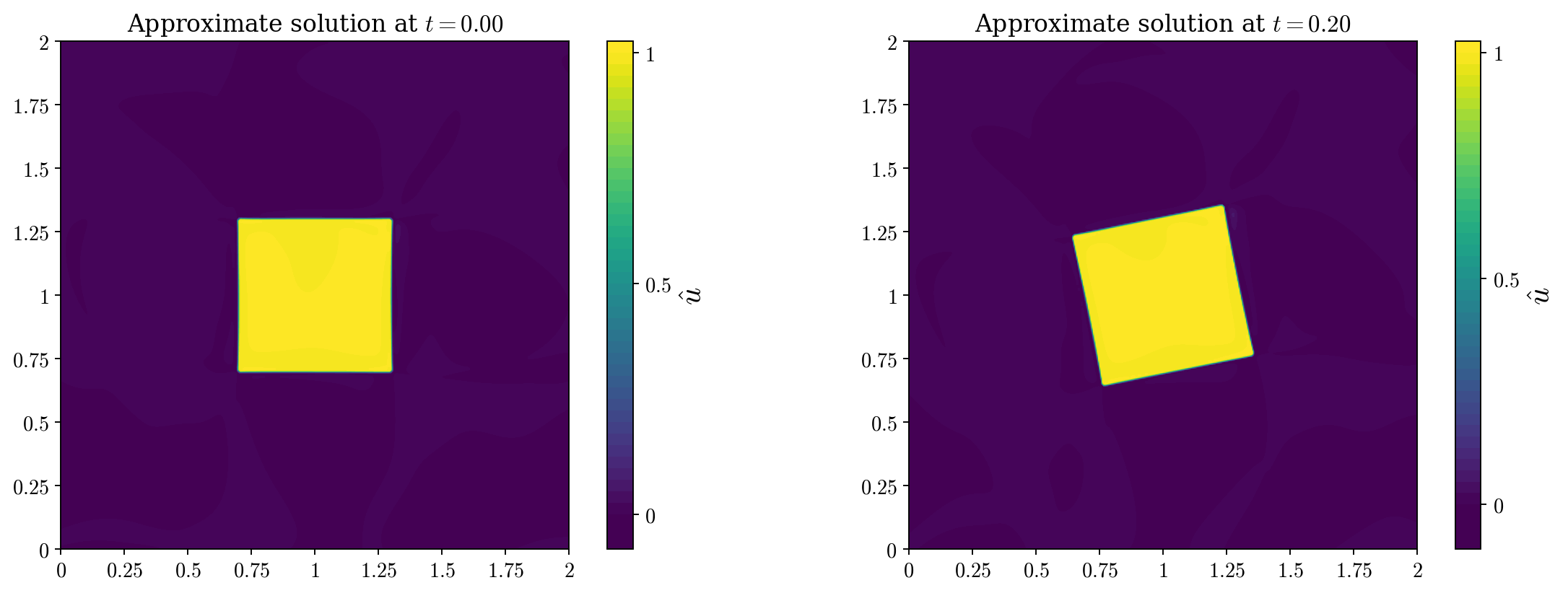}\\[1em]
    \includegraphics[width=\textwidth]{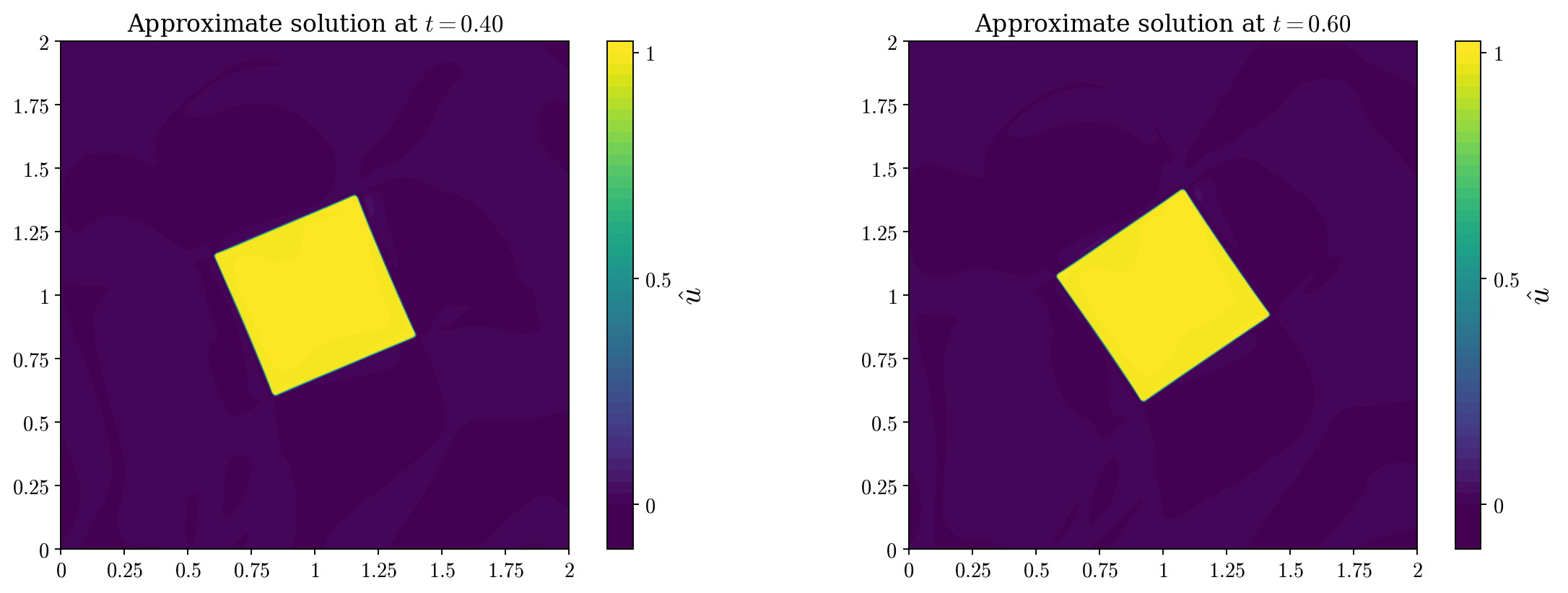}\\[1em]
    \includegraphics[width=\textwidth]{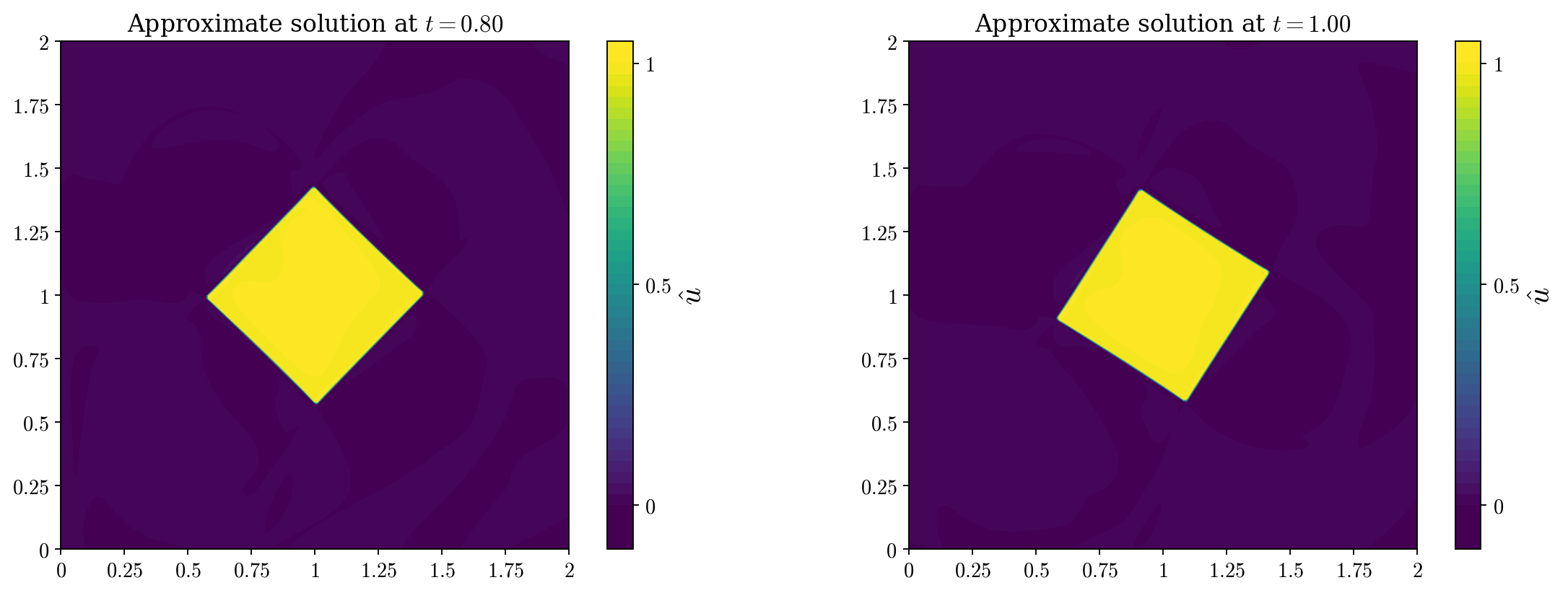}
    \caption{Approximate solution of the 2D rigid‑body rotation at six different times, after applying the corner‑preserving median filter. From top to bottom: \(t = 0.0\) and \(0.2\); \(t = 0.4\) and \(0.6\); \(t = 0.8\) and \(1.0\).}
    \label{fig:all_2d_times}
\end{figure}

Finally, the approximate solution after the corner‑preserving median filter is displayed at six equally spaced times in Fig.~\ref{fig:all_2d_times}. The color scale is fixed between 0 and 1 to emphasis the interface. The sharp square shape is well preserved throughout the rotation.

\subsection{One‑Dimensional Advection of Multiple Blocks}

To isolate and highlight the effects of the median filter, adaptive Fourier weighting, and the two‑stage training of the Fourier mapping, we intentionally employ a simplified setting with insufficient training data and a reduced network capacity. In this experiment, the standard PINN loss based on automatic differentiation is used, and the adaptive sampling strategy is disabled, so that the contribution of each component can be clearly observed. We solve the linear advection equation
\[
u_t + 2u_x = 0, \qquad x\in[0,2],\; t\in[0,1],
\]
with a discontinuous boundary condition that switches from \(0\) to \(1\) at \(t = 0.5\) (\(u(0,t)=0\) for \(t<0.5\) and \(u(0,t)=1\) for \(t\ge 0.5\)) and an initial condition consisting of five rectangular blocks of different heights. The exact solution is the rigid translation of the initial profile together with the propagation of the boundary discontinuity.

Only \(1000\) points are used for the initial condition and \(2000\) interior collocation points — a deliberately low number that causes the raw network output to exhibit strong spurious oscillations. The neural network is a multilayer perceptron with three hidden layers of \(20\) neurons each, and the Fourier feature mapping employs \(D = 40\) sine–cosine pairs. For simplicity, training is performed using only the L‑BFGS optimizer (without the preceding Adam phases).

Since the exact solution contains only jump discontinuities and no point‑like peaks or sharp corners, the corner‑preserving adaptive median filter is unnecessary; a regular median filter suffices. In this test case, the median filter employs a window size \(k = 9\) and is applied to a spatial profile discretized into 1000 uniformly distributed points.

\begin{figure}[H]
    \centering
    \includegraphics[width=0.8\textwidth]{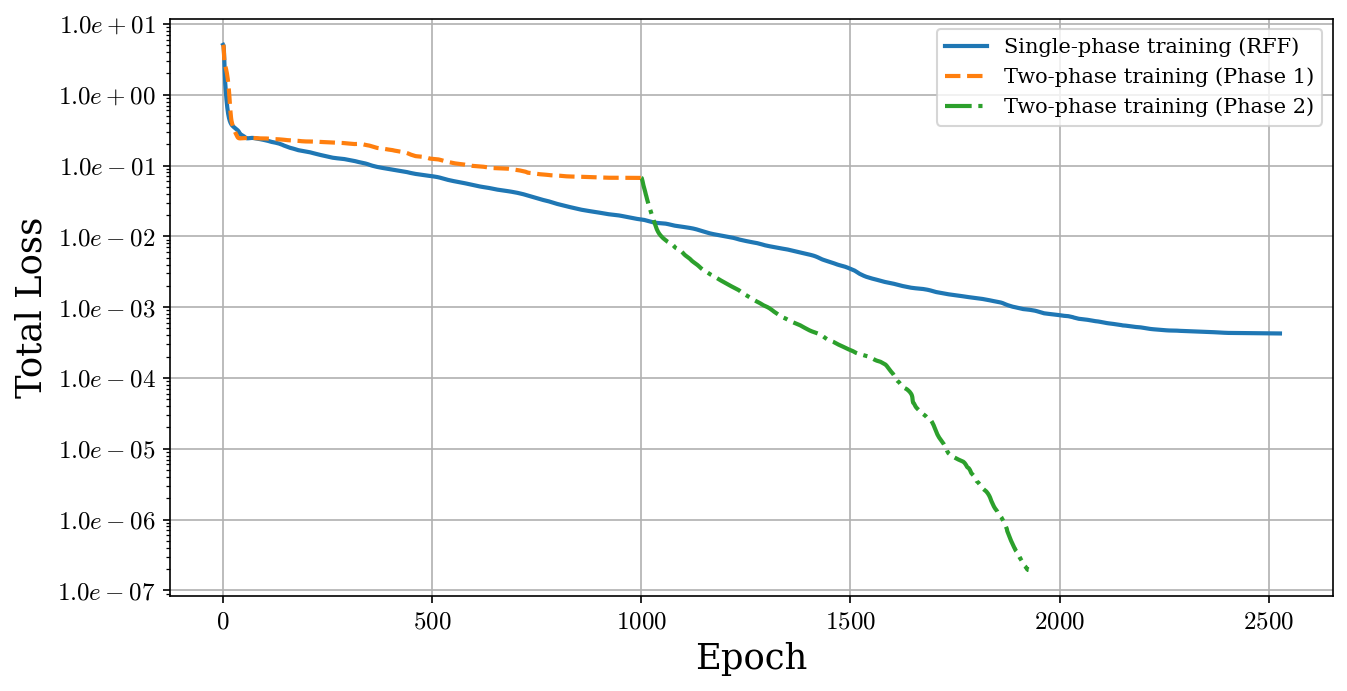}
    \caption{Comparison of trainable versus random Fourier features for the one‑dimensional multiblock problem. Trainable features lead to a noticeable improvement in accuracy.}
    \label{fig:1d_trainable_fourier}
\end{figure}

\medskip
\noindent\textbf{Trainable Fourier features.}
To demonstrate the advantage of trainable Fourier features over the standard random Fourier mapping, we compare the two approaches under identical conditions. Figure~\ref{fig:1d_trainable_fourier} shows the squared \(L^2\) error evolution for a model trained with trainable Fourier features and one trained with fixed random Fourier features (drawn from \(\mathcal{N}(0,\sigma^2)\)). The trainable mapping consistently yields a lower error and faster convergence, confirming that adapting the frequencies to the data is beneficial even in this low‑data regime.

\begin{figure}[H]
    \centering
    \includegraphics[width=0.8\textwidth]{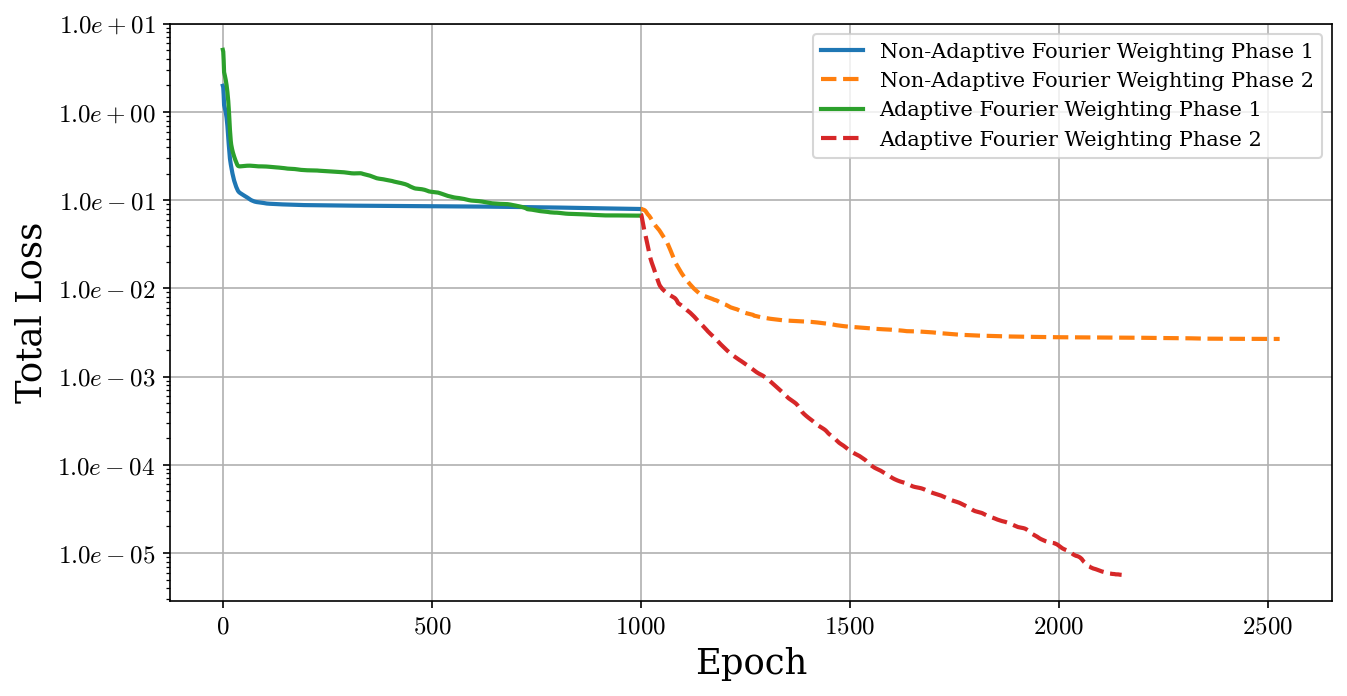}
    \caption{Effect of adaptive Fourier weighting in the two‑stage training. The first stage uses a 20‑fold weight on initial/boundary conditions; the second stage uses gradient‑norm‑based weighting. The adaptive scheme yields faster convergence and a smaller squared \(L^2\) error.}
    \label{fig:1d_adaptive_weighting}
\end{figure}

\noindent\textbf{Adaptive Fourier weighting.}
Next, we examine the effect of adaptive Fourier weighting during the two‑stage training. Figure~\ref{fig:1d_adaptive_weighting} plots the squared \(L^2\) error for the proposed two‑stage scheme with and without the adaptive weighting strategy. In the first stage (Fourier mapping optimization), the loss terms corresponding to the initial and boundary conditions are assigned a weight 20 times larger than that of the PDE residual, forcing the optimizer to prioritize the discontinuous data. In the second stage (network weight training), gradient‑norm‑based weighting is employed. This combined strategy accelerates convergence and results in a lower final error.

\begin{figure}[H]
    \centering
    \begin{subfigure}[b]{0.48\textwidth}
        \centering
        \includegraphics[width=\textwidth]{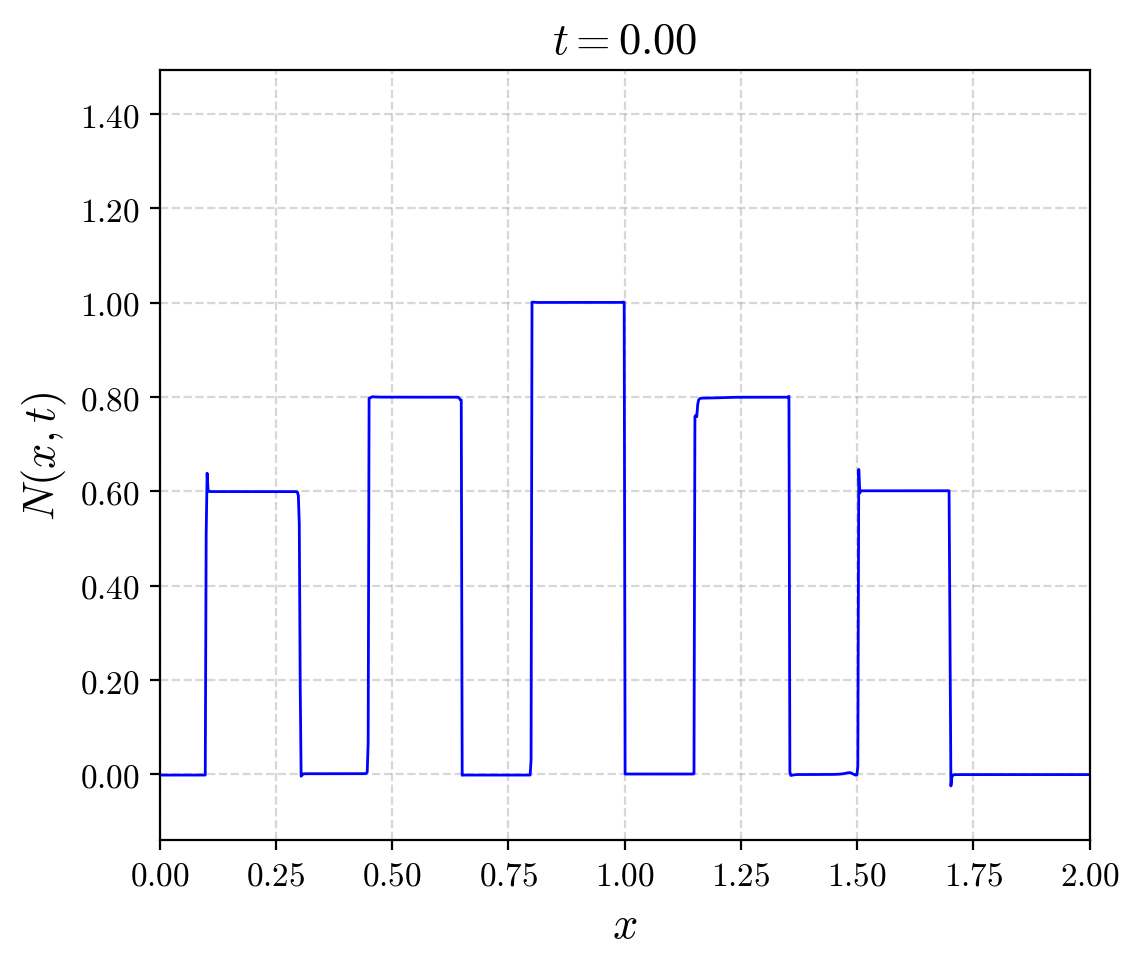}
        \caption{Raw, \(t = 0.2\)}
    \end{subfigure}
    \hfill
    \begin{subfigure}[b]{0.48\textwidth}
        \centering
        \includegraphics[width=\textwidth]{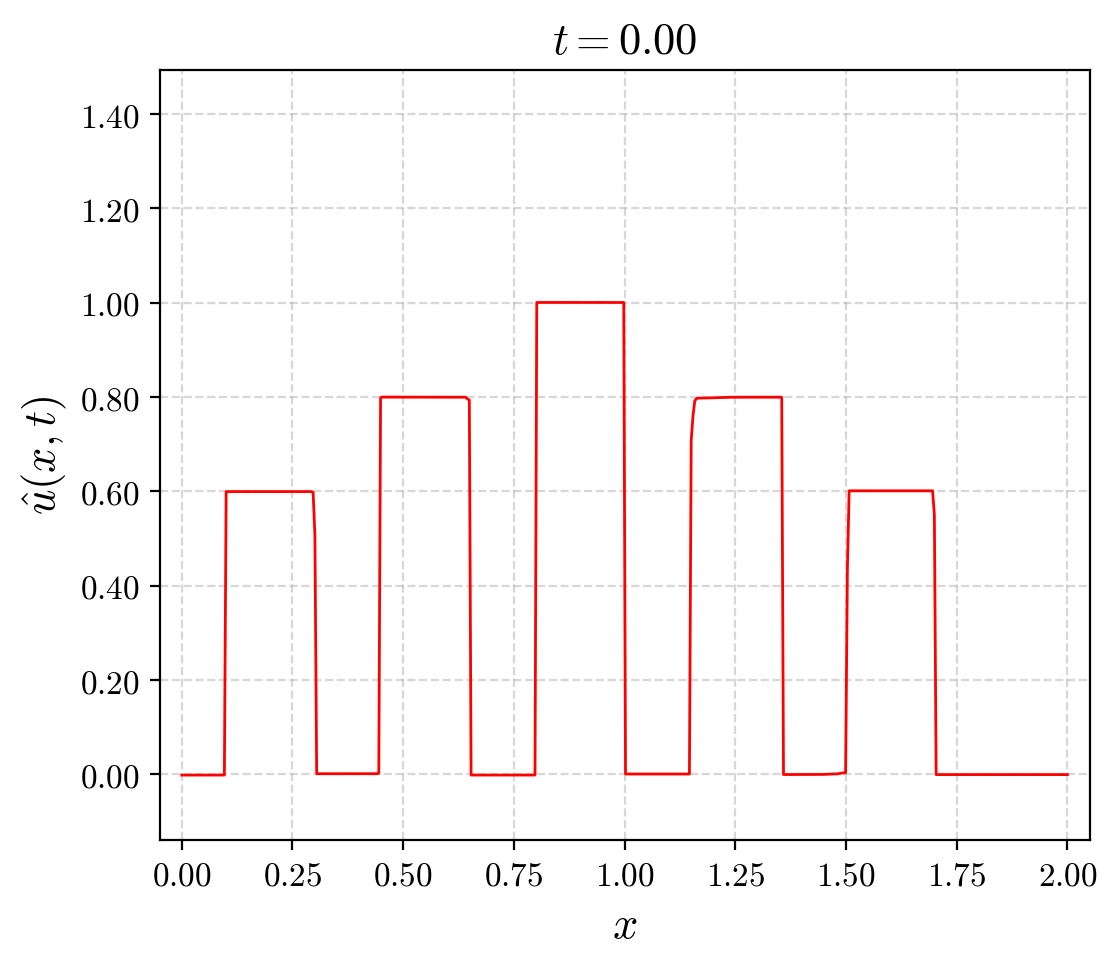}
        \caption{Filtered, \(t = 0.2\)}
    \end{subfigure}

    \vspace{0.5em}

    \begin{subfigure}[b]{0.48\textwidth}
        \centering
        \includegraphics[width=\textwidth]{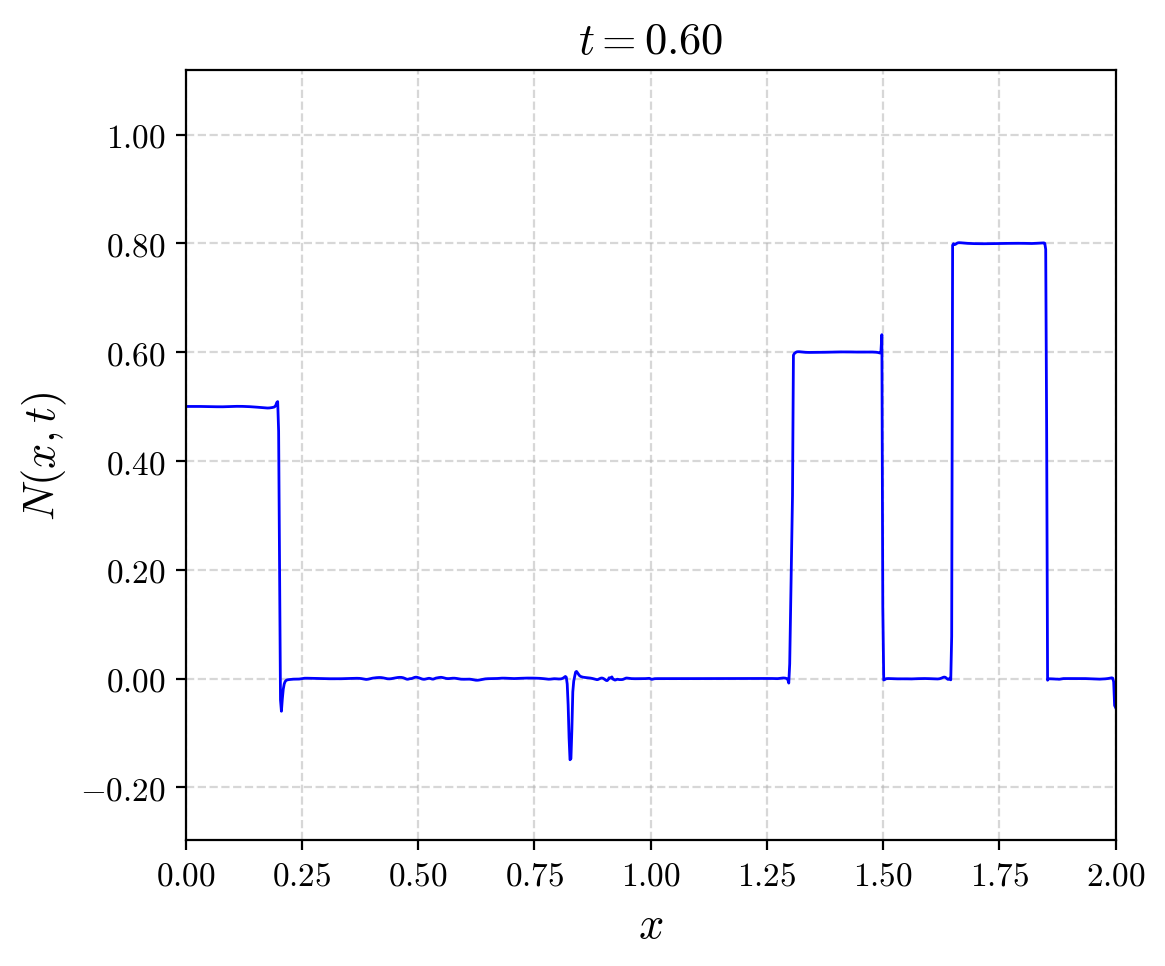}
        \caption{Raw, \(t = 0.5\)}
    \end{subfigure}
    \hfill
    \begin{subfigure}[b]{0.48\textwidth}
        \centering
        \includegraphics[width=\textwidth]{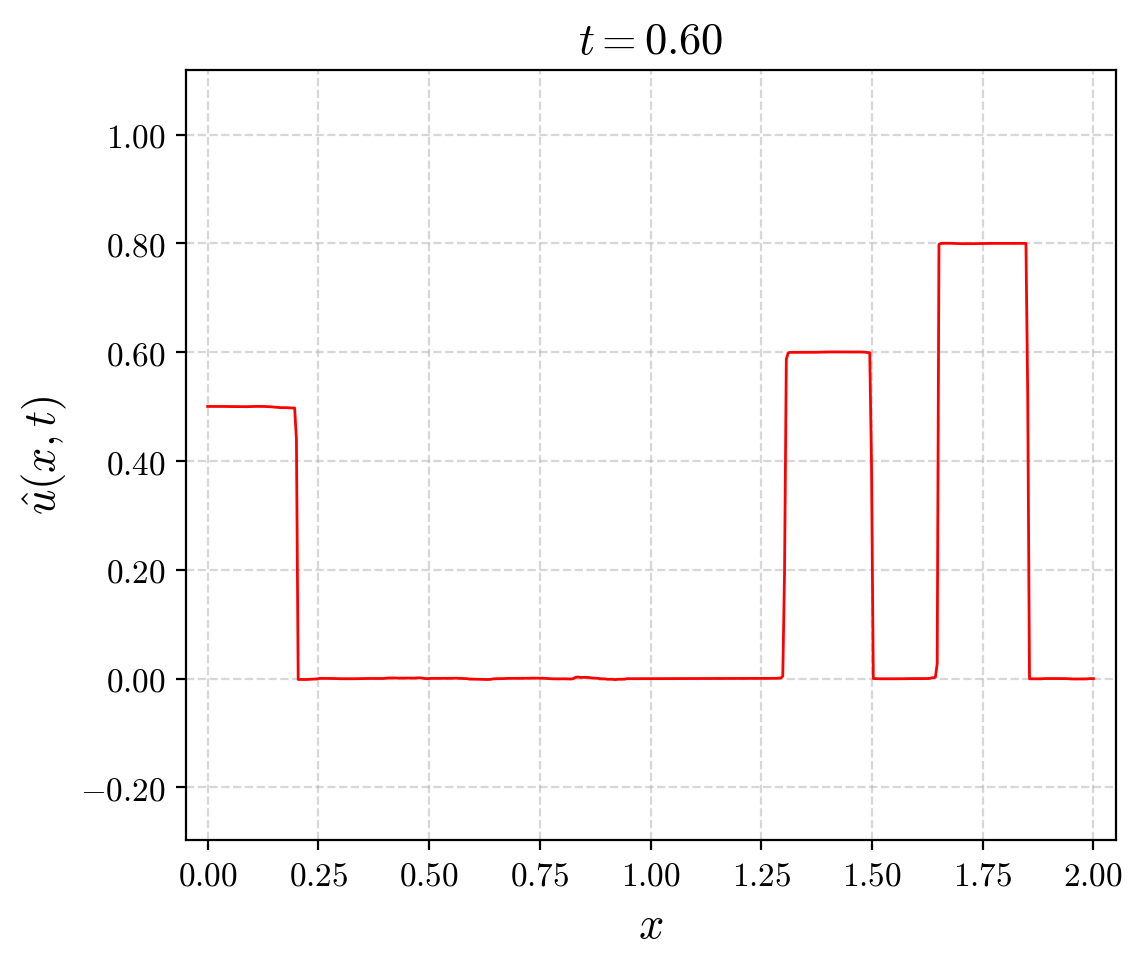}
        \caption{Filtered, \(t = 0.5\)}
    \end{subfigure}

    \vspace{0.5em}

    \begin{subfigure}[b]{0.48\textwidth}
        \centering
        \includegraphics[width=\textwidth]{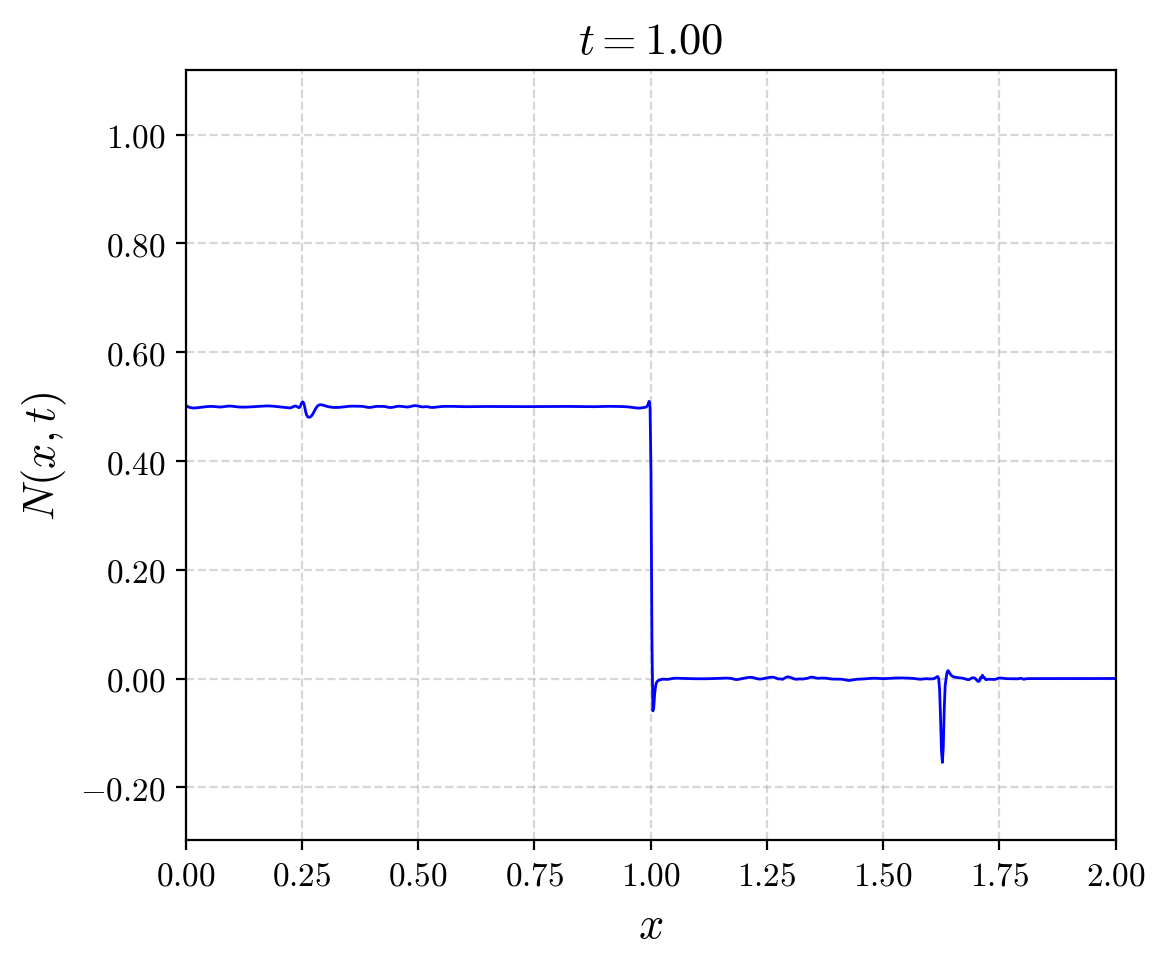}
        \caption{Raw, \(t = 0.8\)}
    \end{subfigure}
    \hfill
    \begin{subfigure}[b]{0.48\textwidth}
        \centering
        \includegraphics[width=\textwidth]{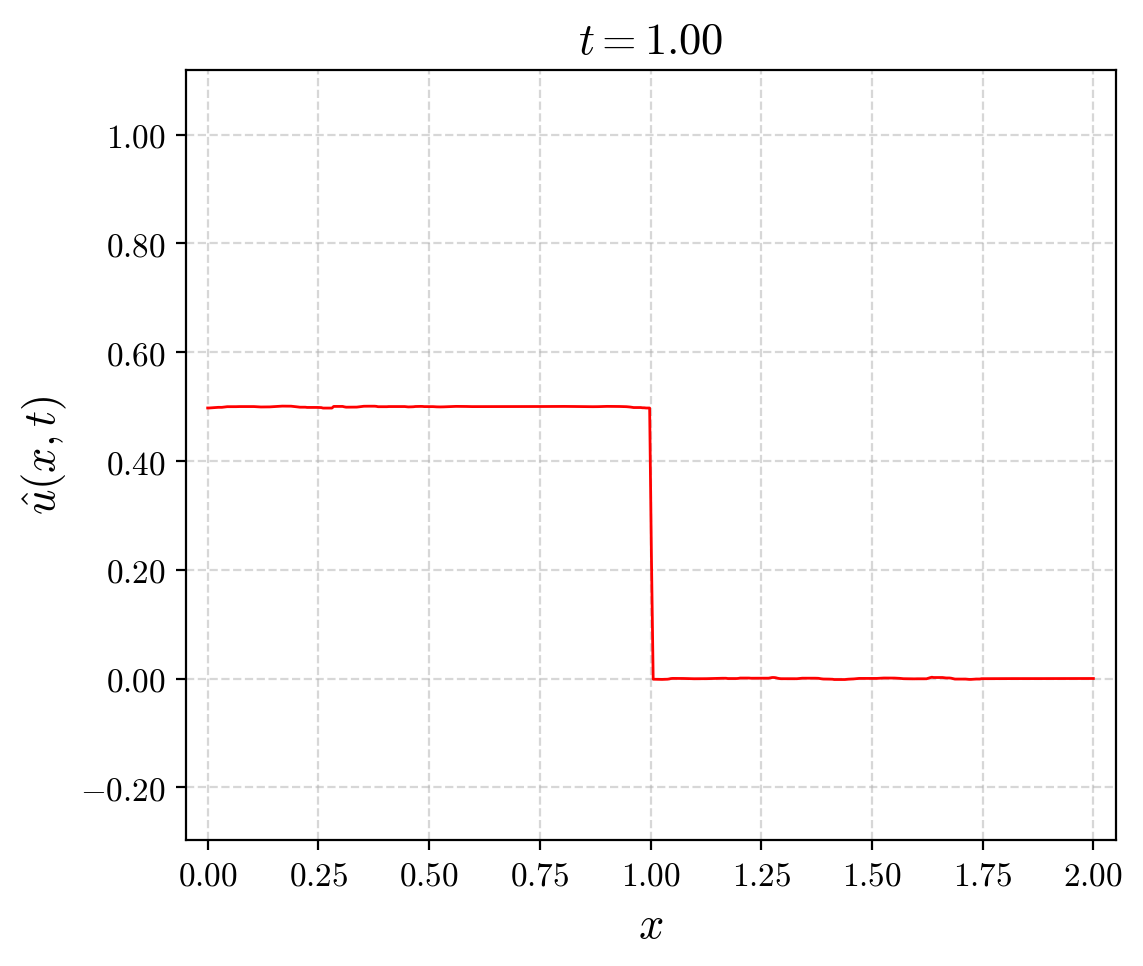}
        \caption{Filtered, \(t = 0.8\)}
    \end{subfigure}

    \caption{Raw network output (left column) and regular median filter result (right column) for the one‑dimensional multiblock problem at three different times. The filter eliminates oscillations while preserving the jumps.}
    \label{fig:1d_filter}
\end{figure}

\medskip
Finally, the combined effect of all components is illustrated in Fig.~\ref{fig:1d_filter}, which shows the raw network prediction and the result after applying the regular median filter (window size \(k = 9\)) at three different times. Despite the extremely limited data, the filter successfully removes the oscillations while keeping the discontinuities sharp, demonstrating the robustness of the proposed approach.

\subsection{Three‑Dimensional Smooth Rotation of a Gaussian Hill}
\label{sec:gaussian_rotation}
To compare the characteristic‑based loss directly with the standard PINN loss in the absence of discontinuities, we consider a fully three‑dimensional smooth problem: rigid‑body rotation of a Gaussian hill about the oblique axis \((1,1,1)\). The advection equation
\[
u_t + (z-y)u_x + (x-z)u_y + (y-x)u_z = 0,
\qquad (x,y,z)\in[-2,2]^3,\; t\in[0,1],
\]
is solved with the initial condition
\[
u(\mathbf{x},0)=\exp\!\big(-20\,\|\mathbf{x}-(1,0,0)\|^2\big).
\]
The domain is chosen large enough that the rapidly decaying bump remains negligibly small near the boundaries throughout the time interval; hence, the homogeneous Dirichlet condition \(u=0\) imposed on all faces is consistent with the exact solution.

Both models use an identical architecture: a multilayer perceptron with four hidden layers of \(20\) neurons, preceded by a Fourier feature mapping with \(D=60\) sine–cosine pairs and \(\sigma=3.0\). Training is performed exclusively with the Adam optimizer (\(3000\) iterations). In addition to the \(N\) interior collocation points listed in Table~\ref{tab:full_comparison}, the training set includes \(2000\) initial condition points and \(300\) boundary points, whose numbers are held fixed across all experiments. For the characteristic‑based loss, the characteristic feet are precomputed with a single fourth‑order Runge–Kutta (RK4) step of size \(h=10^{-2}\), and no automatic differentiation is required for the PDE residual. The standard PINN loss employs automatic differentiation to compute the required spatial and temporal derivatives.

For each value of \(N\), the experiment is repeated with five different random seeds, and the reported metrics are the means over these five runs. The mean squared error (squared \(L^2\) norm) is evaluated on a uniform \(60\times60\times60\) grid at the final time \(t=1\). Table~\ref{tab:full_comparison} shows the per‑iteration time, the total training time, the speed‑up (ratio of standard to characteristic total time), and the squared \(L^2\) error for each data size.

\begin{table}[H]
  \centering
  \setlength{\tabcolsep}{10pt}
  \begin{tabular}{@{}ccccccccc@{}}
    \hline
    \(N\) & \multicolumn{2}{c}{Iter. time (s)} & \multicolumn{2}{c}{Total time (s)} & Speed‑up & \multicolumn{2}{c}{Sq.\ \(L^2\) error} \\
    \cline{2-3} \cline{4-5} \cline{7-8}
    & Char. & Std. & Char. & Std. & (Std./Char.) & Char. & Std. \\
    \hline
    2\,000  & 0.0030 & 0.0066 & 9.2 & 20.3 & 2.20 & 3.50e-4 & 3.56e-4 \\
    8\,000 & 0.0031 &  0.0068 & 9.4 & 20.5 & 2.18 & 1.44e-4 & 1.52e-4 \\
    15\,000 & 0.0031 & 0.0069 & 9.5 & 20.5 & 2.16 & 9.04e-5 & 1.44e-4 \\
    50\,000 & 0.0062 & 0.0102 & 22.9 & 32.8 & 1.43 & 4.90e-5 & 5.72e-5 \\

    \hline
  \end{tabular}
  \caption{Performance metrics for different numbers of interior collocation points \(N\) (all values are averages over five random seeds). ``Char.'' denotes the characteristic‑based loss, ``Std.'' the standard PINN loss (automatic differentiation). The error is the squared \(L^2\) norm evaluated on a \(60\times60\times60\) grid at \(t=1\).}
  \label{tab:full_comparison}
\end{table}

The results indicate that the characteristic‑based loss reduces the training cost without increasing the error, and the speed‑up grows with the data size. This advantage holds even for a completely smooth solution and confirms the effectiveness of replacing automatic differentiation with a characteristic relation.
\subsection{Three‑Dimensional Nonlinear Advection with a Manufactured Smooth Solution}

To demonstrate that the advantage of the characteristic‑based loss extends to nonlinear problems, we consider a fully three‑dimensional nonlinear advection equation with a smooth manufactured solution. The equation
\[
u_t + u\,u_x + u\,u_y + u\,u_z = f(x,y,z,t),
\qquad (x,y,z)\in[-2,2]^3,\; t\in[0,1],
\]
is solved with the initial condition
\[
u(\mathbf{x},0)=\exp\!\big(-20x^{2}-15y^{2}-10z^{2}\big),
\]
which is the manufactured exact solution at \(t=0\). The source term \(f\) is chosen such that the smooth function
\[
u_{\text{exact}}(x,y,z,t)=\exp\!\big(-20(x-0.3t)^{2}-15(y-0.2t)^{2}-10(z-0.1t)^{2}\big)
\]
satisfies the PDE exactly. Substituting this expression into the equation yields
\[
f = u_{\text{exact}}\big[12(x-0.3t)+6(y-0.2t)+2(z-0.1t)\big] - u_{\text{exact}}^{2}\big[40(x-0.3t)+30(y-0.2t)+20(z-0.1t)\big].
\]
Because the Gaussian decays rapidly, the domain is large enough that the solution remains negligibly small on the boundaries; thus the homogeneous Dirichlet condition \(u=0\) on all faces is consistent with the exact solution.

We use the same network architecture and training configuration as in the previous example (Section~\ref{sec:gaussian_rotation}). The characteristic‑based loss computes the characteristic feet on‑the‑fly with a single Euler forward step of size \(h=10^{-2}\), using the current network prediction \(\hat{u}\) as prescribed in~\eqref{eq:char_loss_nonlinear}. No automatic differentiation is required for the PDE residual, whereas the standard PINN loss employs automatic differentiation.

The training set contains \(2000\) initial condition points and \(300\) boundary points, identical to the previous setup, while the number \(N\) of interior collocation points is varied. For each value of \(N\), the experiment is repeated with five different random seeds; the reported metrics are the means over these five runs. The mean squared error (squared \(L^{2}\) norm) is evaluated on a uniform \(60\times60\times60\) grid at the final time \(t=1\).

\begin{table}[H]
  \centering
  \setlength{\tabcolsep}{10pt}
  \begin{tabular}{@{}ccccccccc@{}}
    \hline
    \(N\) & \multicolumn{2}{c}{Iter.\ time (s)} & \multicolumn{2}{c}{Total time (s)} & Speed‑up & \multicolumn{2}{c}{Sq.\ \(L^{2}\) error} \\
    \cline{2-3} \cline{4-5} \cline{7-8}
    & Char. & Std. & Char. & Std. & (Std./Char.) & Char. & Std. \\
    \hline
    2\,000  & 0.0050 & 0.0073 & 15.1 & 22.1 & 1.46 & 1.19e-4 &  4.44e-4 \\
    8\,000  & 0.0053 & 0.0073 &  16.0 &  22.2 & 1.39 &  2.94e-5 & 3.29e-4 \\
    15\,000 & 0.0057 & 0.0079 &  17.4 & 24.0 & 1.38 &  7.11e-6 & 1.72e-5 \\
    50\,000 & 0.0112 & 0.0154 & 33.8 & 46.4 & 1.37 & 3.86e-6 & 2.49e-6 \\
    \hline
  \end{tabular}
  \caption{Performance metrics for the nonlinear manufactured‑solution problem, trained with the Adam optimizer (\(3000\) iterations). ``Char.'' denotes the characteristic‑based loss, ``Std.'' the standard PINN loss (automatic differentiation). All values are averages over five random seeds. The error is the squared \(L^{2}\) norm on a \(60\times60\times60\) grid at \(t=1\).}
  \label{tab:nonlinear_adam}
\end{table}

The results confirm that even in a genuinely nonlinear, smooth setting, the characteristic‑based loss reduces the training cost without sacrificing accuracy. When the L‑BFGS optimizer is employed instead of Adam, the relative speed‑up becomes smaller because the line‑search procedure performs multiple forward evaluations per iteration, which dilutes the savings from avoiding automatic differentiation; nevertheless, the characteristic‑based loss still maintains a slight advantage in total training time. Overall, the elimination of automatic differentiation in the characteristic‑based formulation consistently leads to a noticeable reduction in computational effort.
\section{Conclusion}
\label{sec:conclusion}

In this paper, we introduced a characteristic-based loss function for physics-informed neural networks (PINNs) applied to advection equations. By embedding the characteristic relations directly into the training objective, the proposed loss completely bypasses spatial automatic differentiation, substantially lowering computational cost while retaining the flexibility of the PINN framework. The loss is formulated for both linear and nonlinear advection and provides a general and efficient alternative to the standard PDE residual.

To overcome the specific difficulties associated with discontinuous initial and boundary conditions, we developed a set of complementary techniques: characteristic‑based adaptive sampling that concentrates collocation points near propagating discontinuities, a trainable Fourier feature mapping with two‑stage training and adaptive weighting to reduce spectral bias, and an adaptive spatial median filter that removes spurious oscillations while preserving sharp corners and peaks. These strategies work in concert with the characteristic‑based loss, offering a comprehensive framework for advection problems with sharp fronts.

Extensive numerical experiments on two‑ and three‑dimensional benchmark problems — both linear and nonlinear, smooth and discontinuous — demonstrated that the proposed approach consistently achieves improved or comparable accuracy with lower training cost compared to conventional PINN formulations. The computational advantage of the characteristic‑based loss is particularly pronounced with the Adam optimizer and on GPU hardware, while it remains present with the L‑BFGS optimizer as well. It was also shown that even for completely smooth problems, eliminating automatic differentiation reduces the training time without any sacrifice in accuracy.

These results highlight the effectiveness of exploiting the characteristic structure of the governing equation within the PINN loss function. Combining this loss with tailored sampling and representation techniques for discontinuities significantly improves accuracy and robustness in challenging advection problems.


\end{document}